\documentclass[12pt,a4paper]{article}
\usepackage[cp1251]{inputenc}
\usepackage[english]{babel}
\usepackage{indentfirst}
\usepackage{amsmath}
\usepackage{graphicx}
\usepackage{float}
\usepackage{stmaryrd}
\usepackage{bbm}
\usepackage[all]{xy}
\usepackage{amsfonts}
\usepackage{amssymb}

\usepackage{amscd}
\usepackage{mathrsfs}
\usepackage{tipa}
\usepackage{mathtext,amsmath, amstext, amsgen, amsthm}
\usepackage{amsopn, amsfonts, euscript,amscd,amssymb,latexsym}

\newcommand{\Ind}{\mathop{\rm Ind}\nolimits}
\newcommand{\st}{\mathop{\rm st}\nolimits}
\newcommand{\bk}{\mathop{\rm bk}\nolimits}
\newcommand{\lk}{\mathop{\rm lk}\nolimits}
\newcommand{\bst}{\mathop{\rm bst}\nolimits}
\newcommand{\card}{\mathop{\rm card}\nolimits}

 \newtheorem{theorem}{Theorem}
 \newtheorem{lemma}[theorem]{Lemma}
 \newtheorem{proposition}[theorem]{Proposition}
 \theoremstyle{definition}
 
 \newtheorem{prof}{Proof}

\newtheorem{algorithm}{ALGORITHM}

\begin{document}

\begin{center}
{\bf \large COMPUTING INTERSECTION NUMBERS AND BASES OF COHOMOLOGY GROUPS \\

\vspace{2mm}
FOR TRIANGULATED CLOSED THREE-DIMENSIONAL MANIFOLDS }\\
\end{center}

\begin{center}
{E. I. YAKOVLEV~$^\dag$, V. Y. EPIFANOV~$^\ddag$}\\
\end{center}

{\it $^\dag$, $^\ddag$ Lobachevsky University of Nizhny Novgorod, \\
Department of Information Technology, Mathematics and Mechanics,\\
Gagarina ave., 23, Nizhny Novgorod, 603950, Russia\\
} 
\,\,\,E-mail:{ yei@uic.nnov.ru;\,\,\,vepifanov92@gmail.com} 

\vspace{4mm}
We solve some computational problems for triangulated closed three-dimensional manifolds using groups of simplicial homology and cohomology modulo 2.
Two efficient algorithms for computing the intersection numbers of 1- and 2-dimensional cycles are developed. By means of these algorithms it is
possible to construct a basis of cohomology group from the homology group of two cycles of complementary dimensions.



\section{Introduction}

Computational topology actively develops in last decades and becomes more
and more important in applications (see, for example, Ref. 
1-18). One of the main objectives of this science is the exploration
of methods for calculating topological characteristics of computer models.
The latter are often triangulated topological manifolds.
In algebraic topology, they are called polyhedrons.

In the paper we consider polyhedrons $P$ being closed three-dimensional manifolds.
Our main goal is the elaboration of algorithms for computing
the intersection numbers modulo two of cycles
$x\in Z_m(P)$ and $y\in Z_l(P)$, $m+l=3$.

To achieve this goal is natural to use maps $F_*^{-1}:H_m(P) \to H^l(P)$,
which are the inverse maps for Poincare isomorphisms  $F_*:H^l(P) \to H_m(P)$.

According to the theory inverse map $F_*^{-1}$ is the composition of
isomorphisms $g:H_m(P) \to H_m^*(P)$ and $h:H_m^*(P) \to H^l(P)$,
where $H_m^*(P)$ -- star homology group of a polyhedron $P$.
Unfortunately, there are no known methods for calculating the
map $g$ in general cases, since the existence of the map $g$
is derived from a formal comparison of chain complexes,
which are defined by a group of simplicial homology
$H_m(P)$ and a star homology group $H_m^*(P)$.

Thus, the practical problem of calculating of the
isomorphism $F_*^{-1}:H_m(P) \to H^l(P)$ is still open.
In more detail, it can be reduced to the following problem.
Cycle $x \in Z_m(P)$, consisting of $m$-dimensional simplices
of a simplicial complex $K(P)$ is given. Need to find a set of $l$-dimensional
simplexes of $K(P)$, which barycentric stars form the $m$-dimensional star cycle $x^*$,
homologous to the cycle $x$.

Algorithms 1 and 2 below are the solution of the problem.
In the first algorithm, we assume $m=1$, and $m=2$ for the second one.
The input data is a simplicial structure of the polyhedron $P$ and
the cycle $x \in Z_m(P)$. The output is a co-chain
$J_x:C_1(P) \rightarrow \mathbb{Z}_2$ such
that $\Ind([x],[y])=J_x(y)$ for all $y \in Z_l(P)$.

It is worth noting that similar results for two-dimensional
closed manifolds were obtained in Refs. 7 and 9.
In Ref. 8 cochain $J_x$ is constructed for a given simple cycle
$x \in Z_{n-1}(P)$ on a closed manifold $P$ of arbitrary dimension $n$.
One needs to remark that in the case when cycle $x$ is not simple
the algorithm from Ref. 8 is not applicable. Furthermore, even for $n=3$
and simple cycle $x$ our algorithm essentially differs from that in Ref. 8.

In Theorem 1, the proof of correctness of algorithms is presented, and in Theorems 2 and 3, we calculate their computational complexity.

Algorithms 1 and 2 allow us also to find a basis $[J_{x_1}],\dots,[J_{x_r}]$ of cohomology group $H^l(P)$ with coefficients in $\mathbb{Z}_2$
using a given basis $[x_1],\dots,[x_r]$ of the corresponding homology group $H_m(P)$.

For $m=2$ and $l=1$ cocycles $J_{x_1},\dots,J_{x_r}$ can be used to construct a regular simplicial covering $p:\hat{P}\rightarrow P$
with monodromy group $G \cong H_1(P)$. Covering $p$, in turn, can be used in the problem of minimization of paths and cycles of the manifold $P$
within their homology classes. This approach for $n$-dimensional manifolds was first proposed in Ref. 4 and developed in
Refs. 6 and 8. For the case $n=2$, it is also considered in Refs. 7, 9, 13 and 17.

\section{Preliminaries}

We consider the compact uniform polyhedron $Q$ with the given finite simplicial complex $K(Q)$.
Let us $K'(Q)$ denote a barycentric subdivision of $K(Q)$; $K^m(Q)$ and $K'^m(Q)$ are sets of simplices of dimension $m=0,\dots,n=\dim Q$.
We use simplicial homology groups $H_m(Q)$ and cohomology groups $H^m(Q)$ with coefficients in $\mathbb{Z}_2$,
as well as the corresponding groups of chains and cycles  $C_m(Q)$ and $Z_m(Q)$, co-chains
and co-cycles $C^m(Q)$ and $Z^m(Q)$. Remark that in this case any chain $c \in C_m(Q)$ may be considered as a set of $m$-dimensional simplices,
and as their formal sum.
The union $|c|$ of all simplices from the chain $c$ is called to be its body.

Due to the uniformity of polyhedron $Q$, for any simplex $\sigma \in K^m(Q)$ there is a non-empty set
$\st(\sigma,Q)$ consisting of $n$-dimensional simplices from $K(Q)$ containing $\sigma$.
Obviously, $\st(\sigma,Q) \in C_n(Q)$, and the body $|\st(\sigma,Q)|$ of this chain is the star of the simplex $\sigma$ in $Q$.

Let $v \in K^0(Q)$, $\bk(v,Q)$ be a set of simplices in the boundary $\partial(\st(v,Q))$ containing the vertex $v$,
and $\lk(v,Q)$ be a set of all other simplices of $\partial(\st(v,Q))$.
Then $|\lk(v,Q)|$ is the link of $v$ in polyhedron $Q$, $$\bk(v,Q)=\st(v,\partial(\st(v,Q))) {\text{and}} \partial(\st(v,Q))=\bk(v,Q)+\lk(v,Q).$$

To indicate similar structures in the simplicial complex $ K'(Q)$ we use symbols $\st'(\sigma,Q)$, $\bk'(\sigma,Q)$ and
$\lk'(\sigma,Q)$, respectively. Let us remark that the chain $\st'(\sigma,Q)$ is a part of the barycentric subdivision $\st(\sigma,Q)'$
of chain $\st(\sigma,Q)$ but does not coincide with it. This also holds for $\bk'(\sigma,Q)$.
But chains $\lk'(\sigma,Q)$ and $\lk(\sigma,Q)'$ have no common simplices.

If $ Q $ is a closed manifold, $m \in \{0,\dots,n\}$ and $l=n-m$,
then for each simplex $\sigma \in K^l(Q)$ a chain $\bst(\sigma,Q)$ is
defined, it consists of $m$-dimensional simplices of $K'(Q)$
intersecting with $\sigma$ along its barycenter $\sigma^*$.
The body of this chain is the barycentric star of the simplex $\sigma$ in $Q$.

Let us set
$$
I(\bst(\sigma,Q),\sigma)=1 \in \mathbb{Z}_2, I(\bst(\sigma,Q),\tau)=0 \in \mathbb{Z}_2.
$$
Then for $\tau \in K^l(Q)$, $\tau \ne \sigma$, the formula
$$I(x^*,y)=\sum_{ij} {I(\bst(\sigma_i,Q),\tau_j)},$$
defines the intersection number $I(x^*,y) \in \mathbb{Z}_2$
of the star $m$-dimensional chain $x^*=\sum_i{\bst(\sigma_i,Q)}$
and the simplicial chain $y=\sum_j{\tau_j} \in C_l(Q)$.

If $x^*$ and $y$ are cycles, we can put
$$\Ind([x^*],[y])=I(x^*,y).$$

The last equality correctly defines the bilinear mapping
$\Ind:H_m(Q) \times H_{l}(Q) \rightarrow \mathbb{Z}_2$ which is also called to be the intersection number (Ref. 19, 17.4).

The Poincar\'e isomorphism $F_*:H^{l}(P) \rightarrow H_m(P)$ is induced by
the isomorphism $F:C^l(P) \to C_m(P)$ defined by the formula
\begin{equation}\label{F}
F(J)=\sum_{\sigma \in K^m(P)}J(\sigma)\bst(\sigma,P).
\end{equation}

Thus for any $[J] \in H^l(P)$ and $[y] \in H_l(P)$, we have
\begin{equation}\label {Ind}
\Ind(F_*([J]),[y])=I(F(J),y)=J(y).
\end{equation}

Henceforth $P$ is a 3-dimensional polyhedron being is a closed manifold. In this case,
for each vertex $v \in K^0(P)$, its link $|\lk(v,P)|$ is a triangulated surface
homeomorphic to a sphere $S^2$ (Ref. 20, 2.21). This implies that the star $|\st(v,P)|$ is strongly connected.

\section{Algorithms}

The first algorithm constructs for a given cycle $x \in Z_1(P)$ a co-chain
$J_x:K^2(P) \rightarrow \mathbb{Z}_2$.

\begin{algorithm}
\smallskip

\noindent

\noindent\textbf{Step 1.} For each triangle $t \in K^2(p)$, set $J_x(t):=0$.\\
\textbf{Step 2.} Represent a given cycle as the sum of $x=[v_0,v_1]+\dots+[v_{n-1},v_n]$ of pairwise
different edges, here $v_n=v_0$.\\
\textbf{Step 3.} For each edge $[v_{i-1},v_i]$ from cycle $x$, choose an incident tetrahedron $\sigma_i$. Set $\sigma_0:=\sigma_n$.\\
\textbf{Step 4.} For each number $i=0,\dots,n-1$ in the star of the vertex $v_i$, construct a 3-dimensional path
$\sigma_{i0}+\dots+\sigma_{ik_i}$ from the tetrahedron $\sigma_{i0}=\sigma_i$ to $\sigma_{ik_i}=\sigma_{i+1}$.\\
\textbf{Step 5.} For all triangles $t_{ij}=\sigma_{ij-1} \cap \sigma_{ij}$, where $i = 0,\dots,n-1$
and $j=1,\dots,k_i$ set $J_x(t_{ij}):=J_x(t_{ij})+1$.\\
\end{algorithm}

Further, let $x \in Z_2(P)$ and for any subpolyhedron $Q \subset P$ set
$$
TQ=\{(v,\sigma) \in K^0(Q) \times K^3(P)|\, v \in K^0(\sigma)\}.
$$
The second algorithm constructs mappings $J:TP \rightarrow \mathbb {Z}_2$ and
$J_x:K^1(P) \rightarrow \mathbb{Z}_2$.

\begin{algorithm}
\smallskip

\noindent

\noindent\textbf{Step 1.} For each $(v,\sigma) \in TP$, set $J(v,\sigma):=0$.\\
\textbf{Step 2.} For each vertex $v \in K^0(|x|)$, perform steps 2.1 ~--- 2.3.\\
\textbf{Step 2.1.} Choose an arbitrary tetrahedron $\sigma_0 \in \st(v,P)$ and set
$D:=\bk(v,\sigma_0)$.\\
\textbf{Step 2.2.} For each $\tau \in \bk(v,\sigma_0)$, mark $\sigma_0$ in the list $\partial^{-1}(\tau,P)$.\\
\textbf{Step 2.3.} While $D \ne \emptyset$ perform steps 2.3.1 ~--- 2.3.4.\\
\textbf{Step 2.3.1.} Choose an element $t \in D$ and delete it from the list $D$.\\
\textbf{Step 2.3.2.} If both elements from the list $\partial^{-1}(t,P)$ are marked, then go back to step 2.3.
Otherwise, choose a marked simplex $\sigma_+$ and an unmarked simplex $\sigma$ from the list $\partial^{-1}(t,P)$.\\
\textbf{Step 2.3.3.} If $t \notin x$, set $J(v,\sigma):=J(v,\sigma_+)$. If $t \in x$, set
$J(v,\sigma):=J(v,\sigma_+)+1$.\\
\textbf{Step 2.3.4.} For each $\tau \in \bk(v,\sigma)$, $\tau \ne t$, mark $\sigma$ in $\partial^{-1}(\tau,P)$. If
then in the list $\partial^{-1}(\tau,P)$ a unmarked element exists, add $\tau$ to $D$.\\
\textbf{Step 3.} For each edge $a=[uv] \in K^1(P)$, choose an arbitrary incident simplex
$\sigma \in K^3(P)$ and set $J_x(a):=J(u,\sigma)+J(v,\sigma)$.\\
\end{algorithm}

Let $m=1 $ for the algorithm 1 and $m=2 $ for the algorithm 2, and $l=3-m$.
Then the result of any of these algorithms is the
construction of the co-chain $J_x \in C^{l}(P)$ for the given cycle $x \in Z_m(P)$.
It is supposed co-chain  $J_x: K^{l}(P) \to \mathbb Z_2$ is extended till
the homomorphism $J_x: C_{l}(P) \to \mathbb Z_2$.
\begin{theorem}
For each $m=1,2$, the constructed co-chain $J_x \in C^{l}(P)$ is a cocycle.
The homology class $[x] \in H_m(P)$
and the cohomology class $[J_x] \in H^{l}(P)$ are related by $F_*([J_x])=[x]$,
where $F_*:H^{l}(P) \rightarrow H_m(P)$ is the Poincar\'e isomorphism.
For each cycle $y \in Z_{l}(P)$, the equality $\Ind([x],[y])=J_x(y)$ is valid.
\end{theorem}

\begin{prof} Suppose first $m=1$. Consider in the complex $K'(P)$ the one-dimensional chain \\
\begin{equation}\label{x^*m=1}
x^*=\sum_{i=0}^{n-1}\sum_{j=1}^{k_i}([\sigma_{ij-1}^*t_{ij}^*]+[t_{ij}^*\sigma_{ij}^*])
\end{equation}
and the two-dimensional chain
$$
c=\sum_{i=0}^{n-1}(\sum_{j=1}^{k_i}([v_i \sigma_{ij-1}^* t_{ij}^*]+[v_i t_{ij}^* \sigma_{ij}^*])+[v_{i} v_{i+1} \sigma_i^*]),
$$
where $\sigma_{ij-1}^*$ and $t_{ij}^*$ are barycenters of tetrahedrons $\sigma_{ij-1}$ and triangles
$t_{ij}$ built in steps 4 and 5 of Algorithm 1. Direct calculation shows $\partial c=x+x^*$. Therefore, one has $x^* \in Z_1(P)$ and $[x^*]=[x]$.

In virtue of (\ref{x^*m=1}), we get
\[
x^*=\sum_{i=0}^{n-1} \sum_{j=1}^{k_i}\bst(t_{ij},P).
\]
Moreover, $J_x(t_{ij})=1$ holds for each $i=0,\dots,n-1$, $j=1,\dots,k_i$, and $J_x(t)=0$ for any other
$t \in K^2(P)$. According to (\ref{F}), this implies $F(J_x)=x^*$.

The equality $\partial\circ F=F\circ\delta$ implies $F(\delta J_x)=\partial x^* = 0$.
Since $F:C^l(P) \to C_m(P)$ is isomorphism, it follows that $\delta J_x=0$ and $J_x \in Z^2(P)$.

Due to properties of the Poincar\'e isomorphism $F_*:H^{l}(P) \rightarrow
H_m(P)$ proved above the following equality
$F_*([J_x])=[x^*]=[x]$ is valid. But then, for every cycle $y \in Z_l(P)$, according to
(\ref{Ind}), the following holds
\[\Ind([x],[y])=\Ind([x^*],[y]))=\Ind(F_*([J_x],[y])=J_x(y).\]
This completes the proof of the theorem for $m=1$.
	
For $m=2$, we need in some auxiliary assertions.

\begin{lemma}\label{tau in x criterion}
Suppose that simplices $\sigma, \tilde{\sigma} \in K^3(P)$ contain a common triangle
$\tau=\sigma\cap\tilde{\sigma}$ and $v \in K^0(\tau)$. Then the following is valid:
$\tau \in x$ ($\tau \notin x$), iff $J(v,\sigma)+J(v,\tilde{\sigma})=1$
($J(v,\sigma)=J(v,\tilde{\sigma})$).
\end{lemma}

\begin{prof}
If $v \notin K^0(|x|)$, then one has $\tau \notin x$. However, according to step 1 and conditions in step 2 of the algorithm 2,
we have $J(v,\sigma)=0=J(v,\tilde{\sigma})$.

Suppose further that $v \in K^0(|x|)$. When performing steps 2.1 ~--- 2.3 of the algorithm 2
for any tetrahedron $\theta \in \st(v,P)$, a three-dimensional path
$\mu(\theta)=\sigma_0+\dots+\sigma_q$ in the star $|st(v,P)|$ is implicitly constructed, which the end at $\sigma_q=\theta$,
possessing the following properties:\\
1) for each $i=1,\dots,q$, triangle $t_i=\sigma_{i-1}\cap\sigma_i$ is added to the the list $D$ in step
2.1 or 2.3.4 as the face of the simplex $\sigma_{i-1}$, \\
2) right after selecting the $t_i$ from the list $D$ in step 2.3.1 the tetrahedron $\sigma_{i-1}$ is marked in $\partial^{-1}(t_i,P)$,
and $\sigma_{i}$ is unmarked in the same list. Thus, according to step
2.3.3, one has
$$J(v,\sigma_{i})=J(v,\sigma_{i-1})+I(\bst(t_i,P),x).$$ \\

Denote
$$
\mu^*(\theta)=\sum_{i=1}^{q}\bst(t_i,P).
$$
Then, the previous recurrence relation implies that
\begin{equation}\label{J=I}
J (v, \theta) = I (\mu ^ * (\theta), x).
\end{equation}
Sum $ c = \bst (\tau, P) + \mu ^ * (\sigma) + \mu ^ * (\tilde \sigma) $ is a one-dimensional cycle
lying in the star $ | st (v, P) | $. In this case, $ [c] = 0 $ is in $ H_1 (P) $ and $ I (c, x) = \Ind ([c], [x]) = 0 $.
Consequently, we come to
\begin{equation}\label{I,tau,mu}
I (\bst (\tau, P), x) = I (\mu ^ * (\sigma), x) + I (\mu ^ * (\tilde \sigma), x).
\end{equation}
According to (\ref{J=I}), (\ref{I,tau,mu}), we have
\begin{equation}\label{I,J}
I (\bst (\tau, P), x) = J (v, \sigma) + J (v, \tilde {\sigma}).
\end{equation}
But inclusion $ \tau \in x $ ($ \tau \notin x $) occurs iff $ I (\bst (\tau, P), x) = 1 $ ($ I (\bst (\tau, P), x) = 0 $).
The latter, according to (\ref {I,J}), is equivalent to $ J (v, \sigma) + J (v, \tilde {\sigma}) = 1 $
($ J (v, \sigma) + J (v, \tilde {\sigma}) = 0 $).
\end{prof}

\begin{lemma}\label{correctness}
For any edge $ a = [uv] \in K ^ 1 (P) $ and simplices $ \sigma, \tilde {\sigma} \in K ^ 3 (P) $ incident to $ a $
the equality is valid
\begin{equation}\label{J,u,v}
J(u,\sigma)+J(v,\sigma)=J(u,\tilde{\sigma})+J(v,\tilde{\sigma}).
\end{equation}
\end{lemma}

\begin{prof}
Since $ P $ is a closed manifold then link $ | \lk (a, P) | $ is homeomorphic to the circle (Ref. 20, 2.24).
This means that $ | \st (a, P) | $ is a strongly connected polyhedron. Therefore, it is
sufficient to prove statement of the lemma only for the case where simplices $ \sigma $ and $ \tilde {\sigma} $ contain a common triangle
$ \tau = \sigma \cap \tilde {\sigma} $. But in this situation, by Lemma \ref{tau in x criterion}, when $ \tau \in
x $ we have the equalities $ J (u, \tilde {\sigma}) = J (u, \sigma) + 1 $ and $ J (v, \tilde {\sigma}) = J (v, \sigma) + 1 $
and when $ \tau \notin x $ then equalities $ J (u, \tilde {\sigma}) = J (u, \sigma) $ hold and
$J(v, \tilde {\sigma}) = J (v, \sigma) $. Summed them, we obtain (\ref{J,u,v}) for both cases.
\end{prof}

By Lemma \ref{correctness} the result of the step 3 of the Algorithm 2 does not depend on the choice of simplex $ \sigma $
incident to the edge $ a $.

For $ v \in K ^ 0 (| x |) $, let us $ T ^ + (v) $ be a set of tetrahedrons $ \sigma \in
K ^ 3 (P) $ incident to $ v $ and satisfying the equality $ J (v, \sigma) = 1 $. Denote
\begin{equation}\label{z_0}
z_0 = \sum_ {v \in K ^ 0 (| x |)} \sum _ {\sigma \in T ^ + (v)} {\bk '(v, \sigma)}.
\end{equation}
\begin{equation}\label{z_1}
z_1 = \sum_ {v \in K ^ 0 (| x |)} \sum _ {\sigma \in T ^ + (v)} {\lk '(v, \sigma)}.
\end{equation}

\begin{lemma}\label{t' in z_1 criterion}
For any edge $ a = [uv] \in K ^ 1 (P) $ and the triangle $ t' \in \bst (a, P) $ inclusion $ t' \in z_1 $ is
equivalent to $ J_x (a) = 1 $.
\end{lemma}

\begin{prof}
There is one and only one tetrahedron $ \sigma \in K^3 (P) $ such that $ t' \in K'^{2} (\sigma)$.
Moreover, one has $ t '\in \lk' (u, \sigma) = \lk '(v, \sigma) $.

According to the step 3 of the algorithm 2, the equality $ J_x (a) = 1 $ means that the values of $ J (u, \sigma) $ and $ J (v, \sigma) $ are different,
and therefore simplex $ \sigma $ belongs to only one of two sets $ T ^ + (u) $ and $ T ^ + (v) $. Hence, the triangle $ t '$ belongs to exactly
one term in the right side of (\ref{z_1}), and therefore $ t' \in z_1 $.

If $ J_x (a) = 0 $, then one has $ J (u, \sigma) = J (v, \sigma) $. In this case simplex $ \sigma $ belongs to either both sets $ T ^ + (u) $
and $ T ^ + (v) $ or do not belong to both. Therefore, triangle $t'$ is contained in
an even number of terms in the right side of (\ref {z_1}) and therefore $ t '\notin z_1 $.
\end{prof}

Denote by $K^{1+}(P,x)$ the set of edges $ a \in K^1(P) $ satisfying the equality
$ J_x (a) = 1 $. Let us set
\begin{equation}\label{x^*m=2}
x^*=\sum_{a\in K^{1+}(P,x)}{\bst(a,P)}.
\end{equation}

\begin{lemma}\label{z_0=x' z_1=x^*}
Let $x'$ is the barycentric subdivision of the cycle $x$. Then $z_0 = x'$ and $z_1 = x^*$ are valid.
\end{lemma}

\begin{prof}
Let $P^2$ be 2-dimensional skeleton of $P$. Then for any triangle $t' \in K'^2(P^2)$
there is a unique triangle $t \in K^2(P) $ containing $t'$. Since $P$ is a closed
three-dimensional manifold, then there are exactly two incident to $t$ simplices $ \sigma, \tilde \sigma \in
K^3(P)$. Moreover, one has $t' \in \bk'(v, \sigma)$ and $t' \in \bk'(v, \tilde \sigma)$ where $v$ is a common vertex of
triangles $t'$ and $t$. Since the only one vertex of triangle $t'$ may belong to
the set $K^0(P)$, then for any pair $(w, \hat \sigma) \in TP$, different from $(v, \sigma)$ and
$(v, \tilde \sigma)$, inclusion $t' \in \bk'(w, \hat \sigma) $ is impossible.

According to what has been proved and definition (\ref {z_0}), inclusion $t' \in z_0$ is possible if and only if $ v \in K^0(|x|)$
and only one of two inclusions $\sigma \in T^+(v)$ or $\tilde \sigma \in T^+(v)$ can occur. The last statement is equivalent to
$J(v, \sigma) + J(v, \tilde \sigma) = 1$ which by Lemma \ref {tau in x criterion} may be true if and only if $t \in x $.
Since inclusions $t \in x$ and $t' \in x'$ are equivalent, then the equality $ z_0 = x '$ has been proved.

Consider next the triangle $t' \in K'^2(P) $ which does not belong to the set $K'^2(P^2)$. There is only one simplex $\sigma \in K^3(P)$
containing $t'$ and only one edge $a=[uv] \in K^1(\sigma)$ such that $t' \in \bst(a,P)$. According to (\ref{x^*m=2}) $t' \in
x^*$ may occur
if and only if $ J_x (a) = $ 1. By Lemma \ref{t' in z_1 criterion} the last equality is equivalent to the inclusion $t' \in z_1 $. Thus
we have $z_1=x^*$.
\end{prof}

Now we can prove all statements of the theorem for $m=2$. Indeed, by Lemma \ref{z_0=x' z_1=x^*}, one has $x^*=x' + z_0 + z_1$.
According to (\ref{z_0}) and (\ref{z_1}), we obtain $z_0 + z_1 = \partial {c}$ where
$$
c=\sum_{v\in K^0(|x|)}\sum_{\sigma\in T^+(v)}{\st'(v,\sigma)}.
$$
Hence $x^*$ is a star cycle being homologous to a given cycle $x \in Z_m (P)$.
According to (\ref{x^*m=2}) we have $F(J_x)=x^*$.
But then using exactly the same arguments as in the case $m = 1$, we find that
$ J_x \in Z ^ l (P) $, $ F _ * ([J_x]) = [x ^ *] = [x] $ and $ \Ind ([x], [y]) = J_x (y) $ for every cycle $ y \in Z_l (P)$.
\end{prof}

\section{Computational complexity}

Below the cardinality of a set $ A $ is denoted as $\card A$.
To specify a uniform three-dimensional polyhedron $P$ it is enough to have the lists of vertices $K^0(P)$ and tetrahedrons $K^3(P)$.
Denote $N_3=\card K^3(P)$. Then, the lists of edges $K^1(P)$ and triangles $K^2(P)$ can be constructed for $O(N_3 \log{N_3})$ time,
as well as lists $\partial^{-j}(s,P)$ of incident to them $(i+j)$-dimensional simplices for all simplices $s \in K^i(P)$, $i=0,1,2$,
and numbers $j \in \{1, \dots, 3-i \}$. Since these lists are used in most of computational topology problems,
we will consider them as input of the algorithms 1 and 2.

When evaluating the computational complexity the main parameters will be the numbers $N_3$ and $n_m=\card x$. In the last formula cycle $x$
is considered as a set of $m$-dimensional simplices. We also put $N_i=\card K^i(P)$, $i=0,1,2,3$, and $n_j=\card K^j(|x|)$, $j=0,\dots,m$.

\begin{theorem}
Construction of co-chain $J_x$ using Algorithm 1 has complexity $O(N_3+n_1)$ in the worst case.
\end{theorem}

\begin{prof}
Step 1 can be done in $O(N_2)$ time. Since $N_2 \le 4N_3$, then we have $O(N_2)=O(N_3)$.
Required on the step 2 the representation of a cycle $x$ can be obtained by finding an Eulerian path on the subgraph built on the edges
of that cycle. The running time of the algorithm that finds Eulerian path is $O(n_1)$ (Ref. 21, VI.23). \\

For each edge $a \in x$ as an incident tetrahedron we can select the first element of $\partial^{-2}(a)$ in $O(1)$ time.
Therefore, the total time of step 3 is $O(n_1)$.

Each path constructing in step 4 can be got by using breadth-first search in the abstract graph whose vertices are the simplices of chain $\st(v_i)$,
and edges are pairs of tetrahedrons from $\st(v_i)$ with a common two-dimensional face. In the above graph, each vertex is incident to no more
than three edges. Therefore, breadth-first search can be done in $O(\card \st(v_i))$ time. Since any tetrahedron may be contained in stars of
no more than four different vertices,  we come to inequality $$\sum_{i=0}^{n-1}\card\st(v_i)\le 4N_3.$$
Therefore, the total time of step 4 is not exceed $O(4N_3)=O(N_3)$.

For each path found at the previous step, the numbering can be done in
a time linear in the length of this path. Therefore, the complexity of step 5 is also equal to $O(N_3)$.

Thus, the three steps can be done in $O(N_3)$ time, and two others in $O(n_1)$ time, in the worst case.
As a result, the entire algorithm has complexity $O(N_3+n_1)$.
\end{prof}

\begin{theorem}
Algorithm 2 has the complexity $O(N_3+n_2\log n_2)$ in the worst case.
\end{theorem}

\begin{prof}
At the step 1, for each tetrahedron $\sigma \in K^3(P)$ we consider only 4 pairs $(v,\sigma)$ such that $v \in K^0(\sigma)$.
This requires $O(4N_3)=O(N_3)$ operations.

Since the chain $\st(v,P)=\partial^{-3}(v,P)$ was constructed before, then step 2.1 can be done in $O(1)$ time.

For any tetrahedron $\sigma$ the set $\bk{(v,\sigma)}$ consists of three elements. The set $\partial^{-1}(\tau,P)$
for any triangle $\tau$ consists of two elements. Therefore, step 2.2 can also be realized in $O(1)$ time. The same is true for step 2.3.4.

Complexity of steps 2.3.1 and 2.3.2 is obviously equal to $O(1)$.
At the step 2.3.3,  checking whether the triangle $\tau$ belongs to the list $x$ can be performed in $O(\log n_2)$ time (Ref. 21).
Since step 2 is executed for each vertex $v$ from the cycle $x$ and $\card{K^0(|x|)}=n_0 \le 3n_2$, then the total complexity of step 2
in the worst case is $O(n_2 \log n_2)$.

For each edge, as mentioned before, some simplex incidental to that edge can be found using $O(1)$ operations.
Therefore, for all edges $a$ from $K^1(P)$ the final value of $J_x(a)$ can be obtained in $O(N_1)$ time. Since $N_1 \le 4N_3$, then one has
$O(N_1)=O(N_3)$.

Thus, steps 1 and 3 can be performed in $O(N_3)$ time and step 2 in $O(n_2 \log n_2)$ time. Consequently, the complexity of the whole
algorithm in the worst case is $O(N_3 + n_2 \log n_2)$.
\end{prof}

\section{Applications}

According to Theorems 1 -- 3, algorithms 1 and 2 allow us for any cycles $x \in Z_m(P)$ and $y \in Z_l(P)$, $m=1,2$, $l=3-m$,
to compute effectively the intersection number using the formula $\Ind([x],[y])=J_x(y)$.
Theorem 1 also implies

\begin{theorem}
If $[x_1],\dots,[x_r]$ is a basis of the homology group $H_m(P)$, $m=1,2$, and $J_{x_1},\dots,J_{x_r}$ are cocycles constructed from
cycles $x_1,\dots,x_r$ using algorithms 1 and 2, then the cohomology classes $[J_{x_1}],\dots,[J_{x_r}]$ form a basis of $H^{l}(P)$, $l=3-m$.
\end{theorem}

Thus, algorithms 1 and 2 can be used to calculate the bases of the cohomology groups of an arbitrary closed three-dimensional manifold.

The collection $J=(J^1,\dots,J^r)$, where $J^i=J_{x_i}$ for all $i=1,\dots,r$, generates a homomorphism $J:C_l(P) \rightarrow \mathbb{Z}_2^r$.
Vector $J(y)$ is naturally called to be the index of the chain $y \in C_l(P)$ with respect to the basis $[x_1],\dots,[x_r]$.

Setting $J([y])=J(y)$ for all $y \in Z_l(P)$, we also obtain a homomorphism $J:H_l(P) \rightarrow \mathbb{Z}_2^r$.

According to the Poincar\'e duality $H_m(P) \cong H_l(P)$, there exists a basis $[y_1],\dots,[y_r]$ of the homology group $H_l(P)$,
for which $\Ind([x_i],[y_j])=\delta_{ij} $, where $\delta_{ij}=1$ for $i=j$ and $\delta_{ij}=0$ for $i \ne j$.
In this case, one has $J([y_j])=e_j$ for all $j=1,\dots,r$, where $e_1,\dots,e_r$ is the standard basis of the group $\mathbb{Z}_2^r$.
Hence

\begin{proposition}
If $y,z \in C_l(P)$ and $\partial{y}=\partial{z}$, then the chains $y$ and $z$ are homologous iff $J(y)=J(z)$.
\end{proposition}

Thus, the algorithms 1 and 2 can be used to solve the problem on if two chains and cycles are homologous.

Suppose further that $m=2$ and $l=1$. Denote $\hat{V}=K^0(P)\times\mathbb{Z}_2^r$.
A collection $\{\hat{v}_i=(v_i,k_i)\in\hat{V}|i=0,\dots,j\}$, $0 \le j \le 3$, will be included in the list $\hat{K}$ when $j=0$, as well as
for $j>0$ and the following conditions: $[v_0 \dots v_j] \in K^j(P)$ and $k_i=k_0+J([v_0v_i])$ for all $i=1,\dots,j$.
The pair $\hat{S}=(\hat{V},\hat{K})$ constructed is an abstract simplicial scheme.
If $\hat{P}$ is any of its realizations, then we have $K^0(\hat{P})=\hat{V}$ and $K(\hat{P})=\hat{K}$.

The natural projection $p:K^0(P) \times \mathbb{Z}_2^r \rightarrow K^0(P)$ is extended uniquely up to a simplicial mapping
$p:\hat{P} \rightarrow P$. According to Ref. 8, the following holds

\begin{proposition}
Mapping $p:\hat{P} \rightarrow P$ is a simplicial regular covering with a monodromy group $G \cong H_1(P)$.
\end{proposition}

Let $y=[v_0v_1]+\dots+[v_{q-1}v_q]$ be a path in $P$. Then for any vector $k_0 \in \mathbb{Z}_2^r$ the only path $\hat{y}$ in $\hat{P}$
which starts at $\hat{v}_0=(v_0,k_0)$ and covers $y$, has the form $\hat{y}=[\hat{v}_0\hat{v}_1]+\dots+[\hat{v}_{q-1}\hat{v}_q]$,
here for each vertex $\hat{v}_i=(v_i,k_i)$ the relation $k_i=k_0+J([v_0v_1]+\dots+[v_{i-1}v_i])$, $i=1, \dots, q$, holds.
In particular, one has $k_q=k_0+J(y)$.

Consider another path $z$ in $P$ from the vertex $v_0$ to $v_q$ and its covering path $\hat{z}$ in $\hat{P}$ which starts at $\hat{v}_0$.
Then the last vertices of the paths $\hat{z}$ and $\hat{y}$ coincide if and only if $J(z)=J(y)$.
By Proposition 1, this is equivalent to the homology of the paths $z$ and $y$.

If $L:C_1(P) \rightarrow \mathbb{R}$ is a non-negative weight function, then the equation $\hat{L}=L \circ p$ determines a similar weight
function $\hat{L}:C_1(\hat{P}) \rightarrow \mathbb{R}$.
Due to that proved above the following holds

\begin{proposition}
Path $z$ in $P$ from the vertex $v_0$ to $v_q$ has the minimal weight $L(z)$ among all paths of the polyhedron $P$, homologous to $y$,
iff the covering path $\hat{z}$ has minimal weight of $\hat{L}(\hat{z})$ among all paths of the polyhedron $\hat{P}$,
going from $\hat{v}_0$ to $\hat{v}_q$.
\end{proposition}

Thus, the problem on the conditional minimum in $P$ is equivalent to the problem of an absolute minimum on the covering polyhedron $\hat{P}$.
In the case of small values of rank $r$ of homology groups $H_1(P)$ and $H_2(P)$, this leads to an algorithm for minimizing the path $y$
in its homology class. In this case, the construction of the polyhedron $\hat{P}$ is actually not required. The algorithm uses only
its above described simplicial scheme. Recall that the algorithm in Ref. 8 requires of simplicity of the initial cycles
$x_1,\dots,x_r$ that form a basis of $H_1(P)$. In this paper, this restriction is overcome.

Unfortunately, the situation is worse for the general case, since, by the construction, $\card \hat{V}=\card{K}^0(P) 2^r$.
Therefore, all known algorithms for the solving the reduced problem for $\hat{P}$ will have an exponential
in the parameter $r$ complexity.

\section{Acknowledgement}

Authors thank the Ministry of Science and Education of Russia for a
support (project 1410, the State Target).

\end{document}